\documentclass{article}
\usepackage{amsmath,amsthm,amssymb,graphics,array}

\theoremstyle{plain}
\newtheorem{theorem}{Theorem}
\newtheorem{lemma}{Lemma}
\newtheorem{proposition}{Proposition}

\newtheorem{corollary}{Corollary}

\theoremstyle{definition}

\newcolumntype{L}{>{$}l<{$}}
\newcolumntype{C}{>{$}c<{$}}

\title{Tamari lattices, forests and Thompson monoids}
\date{\small Mathematics Subject Classification: 06A07, 20F55, 05C05 }

\author{Zoran \v{S}uni\'{c}\\
\small Department of Mathematics\\[-0.8ex]
\small Texas A\&M University, MS-3368\\[-0.8ex]
\small College Station, TX 77843-3368, USA\\[-0.8ex]
}

\begin{document}
\maketitle

\begin{abstract}
A connection relating Tamari lattices on symmetric groups regarded
as lattices under the weak Bruhat order to the positive monoid $P$
of Thompson group $F$ is presented. Tamari congruence classes
correspond to classes of equivalent elements in $P$. The two well
known normal forms in $P$ correspond to endpoints of intervals in
the weak Bruhat order that determine the Tamari classes. In the
monoid $P$ these correspond to lexicographically largest and
lexicographically smallest form, while on the level of permutations
they correspond to $132$-avoiding and $231$-avoiding permutations.

Forests appear naturally in both contexts as they are used to model
both permutations and elements of the Thompson monoid.

The connection is then extended to Tamari orders on partitions of
$((k-1)n+2)$-gons into $(k+1)$-gons and Thompson monoids $P_k$, $k
\geq 2$.
\end{abstract}

\section{Introduction}

The purpose of this note is to present a connection between the
positive Thompson monoid (of type $F$) and Tamari lattices (of type
$A$).

The fact that Thompson groups and monoids are related to trees is
certainly well known and established among the people interested in
these groups. Ever since~\cite{brown:finiteness} trees are heavily
used as helpful tools in representing elements in order to aid both
calculations and conceptual understanding. On the other hand,
researchers in combinatorics have usually, with rare exceptions,
heard very little about Thompson groups, and even when they have it
is mostly in the context of providing examples of infinite simple
groups.

We start with some very well known and understood classes of objects
in combinatorics that are related to Tamari lattices on symmetric
groups and then naturally arrive at Thompson monoids, which
essentially capture all instances of these combinatorial objects
along with their inter-relations.

The connection in question relates Tamari lattices on finite
symmetric groups (Tamari lattices of type $A$) to the positive
Thompson monoid
\begin{equation}\label{presentation-p2}
 P_2 = Mon \langle \ x_0,x_1,x_2,\dots \mid x_i x_j = x_{j+1}x_i, \text{ for }i<j \ \rangle.
\end{equation}

The connection is obtained in a natural way as follows. First some
well known connections between permutations, inversion sequences and
linearized labeled binary rooted trees are recalled. The simple
observation that concatenation is closed in the set $X_\infty$ of
inversion sequences leads to a definition of a graded product on the
set of all finite permutations $S_\infty$. The corresponding product
on the set of linearized labeled binary rooted trees $T_\infty$ is
just stacking of trees. At this stage we have three isomorphic
monoids $X_\infty$, $S_\infty$ and $T_\infty$. Tamari congruence on
$T_\infty$ is the congruence obtained when trees that have the same
shape but different linearization are identified. This leads to a
corresponding congruence on the level of permutations and also on
the level of inversion sequences and we get three monoids
$T_\infty/{\sim}$, $S_\infty/{\sim}$ and $X_\infty/{\sim}$. It turns
out that these three monoids are free. We then extend our
considerations to the set of all sequences of non-negative integers
$X^*$, all linearized labeled binary rooted forests $T^*$ and the
corresponding set of $*$-permutations $S^*$. We extend the notion of
Tamari congruence and identify two forests of the same shape
regardless of the linearization. The corresponding factor monoids
$T^*/{\sim} \cong S^*/{\sim} \cong X^*/{\sim}$ are isomorphic to the
Thompson monoid $P_2$.

After going through the details of the connection between Tamari
lattices on symmetric groups and Thompson monoid $P_2$ in
Section~\ref{delinearization}, a similar connection is established
between Thompson monoids $P_k$, $k \geq 2$, and Tamari orders (they
do not form lattices for $k \geq 3$) corresponding to partitions of
$((k-1)n+2)$-gons into $(k+1)$-gons in Section~\ref{partitions}.

\section{Some basic facts about Tamari lattices}\label{fact-tamari}

Tamari lattices of type $A$ are homomorphic images of the weak
Bruhat order lattices over finite Coxeter groups of type $A$, i.e.,
finite symmetric groups. Recall that the (left) weak Bruhat lattice
on $S_n$ as Coxeter group of type $A_{n-1}$ is just the (left)
Cayley graph of $S_n$ with respect to the standard generating set of
reflections $\{ (12),(23),\dots,(n-1 \ n)\}$ ordered by declaring
that $\sigma \preceq \rho$ if there exists a geodesic path from $1$
to $\rho$ that passes through $\sigma$. Alternatively, we may say
that $\sigma$ is covered by $\rho$ if $\rho = (i \ i+1) \circ
\sigma$, for some standard reflection $(i \ i+1)$, and the length of
$\rho$ (in terms of the standard reflections) is larger than the
length of $\sigma$. Then define the weak Bruhat order as closure of
this cover relation.

For a fixed $n$, there are many ways of thinking of the congruence
classes on $S_n$ defining the corresponding Tamari lattice $L_n$. We
recall some of them here, along with some additional notions.

\begin{figure}[!ht]
\begin{center}
\includegraphics{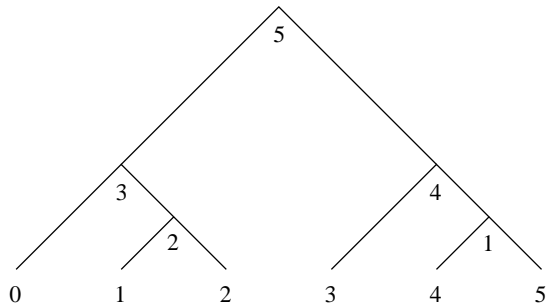}
\caption{A linearized tree}\label{tree-lin}
\end{center}
\end{figure}

We consider linearized labeled binary rooted trees on $n$ interior
vertices. When $n=0$ there is only one such tree and it has a single
vertex which is simultaneously the root and a leaf labeled by 0. If
$n \geq 1$ the root of such tree has degree 2 and the other $n-1$
interior vertices have degree 3. The $n+1$ leafs are labeled
bijectively by $0,1,\dots,n$. In addition, the interior vertices are
labeled bijectively by $1,\dots,n$ in such a way that the labels on
the paths from the root to the leafs are decreasing (this is the
linearization part of the tree - we can use it to extend the partial
order on the interior vertices induced by the tree structure to a
linear order). We depict such trees as in Figure~\ref{tree-lin}.
A rooted binary tree often admits more than one linearization and
some standard choices are well established. We mention here two such
choices (which are relevant to our discussion). The post-order
linearization labels the interior vertices by $1,2,\dots,n$ exactly
in the order they are visited by using the left-right-root rule. The
inverse post-order (or the right-left post order) labels the
interior vertices in the order they are visited by using the
right-left-root rule. The in-order labeling (using the
left-root-right) does not necessarily produce a proper
linearization. In the rest of the text we often say linearized tree
and tree when we mean linearized labeled binary rooted tree and
labeled binary rooted tree (the latter lack linearization labels,
i.e.~they lack labels on the interior vertices).

We recall the interpretation of Tamari lattice $L_n$ as given by
Huang and Tamari in~\cite{huang-t:lattice}. It is defined by the
product order on the set of integer sequences $e_1 \dots e_n$ such
that $1 \leq e_i \leq n$, for all $i$, and the condition that
whenever $i<j$ and $j \leq e_i$ then $e_j \leq e_i$. While this is
not directly apparent in~\cite{huang-t:lattice}, one can easily
interpret these sequences as encodings of labeled binary rooted
trees as follows. Let $t$ be a tree with $n$ interior vertices. For
$i=1,\dots,n$, let $e_i$ be the largest leaf label of the subtree of
$t$ hanging below the interior vertex $i$ in the in-order labeling.
For example, the encoding of the tree in Figure~\ref{tree-lin} is
$22555$.

The Tamari lattice $L_n$ is defined by Bj{\"o}rner and Wachs
in~\cite{bjorner-w:shellable2} as the product order on the set of
integer sequences $r_1 \dots r_n$ satisfying $0 \leq r_i \leq n-i$,
$i=1,\dots,n$, and $r_{k+i} \leq r_k - i$, for $k=1,\dots,n-2$ and
$i=1,\dots,r_k$. A correspondence is established between
permutations and labeled trees (the labeling on the interior
vertices does not respect the partial order imposed by the tree, but
it is related to it in a different way). All permutations are
encoded by integer sequences of the above type as follows. Given a
permutation $\sigma$, for each $i$, count the number of consecutive
terms in $\sigma^{-1}$ following $\sigma^{-1}(i)$ that are smaller
than $\sigma^{-1}(i)$. For example, for $\sigma = 52143$ we have
$\sigma^{-1} = 32541$ and the encoding sequence is $10210$. In the
corresponding tree this sequence records, for each $i$, the number
of interior vertices in the right subtree below the vertex visited
at position $i$ using the in-order. For example, the tree in
Figure~\ref{tree-lin} (ignore the labels on interior vertices) is
encoded by the sequence $10210$. This same tree encodes the
permutation $\sigma=52142$. The Tamari classes are then classes of
permutations encoded by the same integer sequence. The top
permutation in each of these classes is 312-avoiding and the bottom
one is 132-avoiding.

A correspondence between linearized binary trees on $n$ interior
vertices and permutations in $S_n$ is given by Loday and Ronco
in~\cite{loday-r:hopftrees}. The Tamari classes correspond to
classes of permutations that are associated to the same tree
(ignoring the linearization). This is exactly the way in which we
will think of Tamari congruence classes on $S_n$.

There is a way to define triangulations of a $(n+2)$-gon
corresponding to permutations in $S_n$. The Tamari congruence
classes then consists of permutations that produce the same
triangulation (see~\cite{edelman-r:stasheff}). The partial order on
triangulations inherited from the weak Bruhat order is actually
rather natural and can be expressed in its own right, with no
reference to the weak Bruhat order (the cover relation expresses a
local change in the triangulation due to a single ``diagonal
flip''). This is precisely defined in a more general setting in
Section~\ref{partitions}.

Purely in terms of the weak Bruhat order one can define the Tamari
congruence on $S_n$ as the coarsest congruence ${\sim}$ that
collapses the edges
\[
\begin{matrix} (i \ i+1) (i+1 \ i+2) \\ | \\ | \\ (i+1 \ i+2)
\end{matrix},
\]
for $i=1, \dots, n-2$, in the weak Bruhat order on $S_n$. This and
many other lattices on $S_n$ are described in this manner by
N.~Reading in~\cite{reading:cambrian}. The collapsing of edges is
encoded in the Coxeter diagram $A_{n-1}$ by directing the edges. The
Tamari congruence corresponds to orienting all the edges the same
way, as in
\[ (12) \longrightarrow (23) \longrightarrow \dots \longrightarrow (n-1 \ n). \]

There are many other ways to arrive at an ordered lattice isomorphic
to the Tamari lattice corresponding to $S_n$, with or without
referring to permutations. The author has stumbled upon yet another
way in~\cite{sunik:sds} in which fixed points of a certain
endomorphism of an infinite rooted tree are studied. Each vertex
stabilizes after finitely many applications of the endomorphism and
reaches a, so called, self-describing sequence. Each class of points
at level $n$ that eventually stabilizes to the same self-describing
sequence corresponds in a rather natural way (through site inversion
counting) to a congruence class in the Tamari lattice on $S_n$.

Note that there are certainly different congruences on $S_n$
producing the same lattice quotient and thus deserving of the title
Tamari congruence. The point is that there are always some choices
involved and there is often more than one natural choice. One could
work with the right Bruhat order instead of the left one, or define
slightly different way to associate triangulations to permutations,
or apply some obvious automorphisms to the weak Bruhat order
lattice, and so on. There is just too much symmetry involved to
claim any canonical choices (in our discussion so far we already
mentioned a few different choices existing in the literature).

We now fix a particular Tamari congruence on $S_n$. The congruence
will be denoted by ${\sim}$. It is the one we already defined above
in terms of collapsing edges in the weak Bruhat order. We will make
our definition of a triangulation corresponding to a permutation
consistent with this choice. We will also make all our subsequent
choices in accordance to this choice. This makes all the connections
we display possible at the price of not always choosing the most
standard way of representing some objects. It is all matter of left
versus right, bottom versus top, etc., and it seems a standard
choice in one aspect leads to non-standard choices in another
aspect, so some degree of ``oddness'' is unavoidable.

Lattice congruence classes in finite lattices always form intervals,
so the Tamari congruence classes are intervals in the weak Bruhat
order. In agreement with our choices Tamari congruences correspond
to intervals in the weak Bruhat order in which the upper bound is a
$132$-avoiding permutation and the lower bound is a $231$-avoiding
permutation. We will later say more on pattern avoidance.

\section{Some basic facts about Thompson monoid $P_2$}

Here we list some basic facts on Thompson group $F_2$ and the
related positive monoid $P_2$. Other than the presentation, already
given in~\eqref{presentation-p2}, nothing in this section is needed
to follow the text.

The monoid $P_2$ satisfies the Ore condition and embeds into its
group of left fractions $F_2=P_2^{-1}P_2$. The group $F_2$ is the
celebrated Thompson group $F$, given by the group presentation
\[
 F_2 = Gp \langle \ x_0,x_1,x_2,\dots \mid x_i x_j = x_{j+1}x_i, \text{ for }i<j \ \rangle,
\]
which looks exactly the same as the monoid presentation for $P_2$.
The monoid $P_2$ is just the positive submonoid of $F_2$, i.e.~ the
submonoid generated by the set $X=\{x_0,x_1,x_2,\dots\}$.

The element $x_n$, $n \geq 1$, in this presentation can be written
it terms of $x_0$ and $x_1$ as $x_n = x_0^{n-1}x_1x_0^{-(n-1)}$.
Thus $F_2$ is a finitely generated group. It is also finitely
presented (with only two relations), but is it often more convenient
to use the above infinite presentation.

The group $F_2$ has many fascinating properties and has been
studied and re-discovered many times in the last 40 years. It has
been a steady source of highly non-trivial and important examples
or counter-examples, especially in topology.

The group $F_2$ is infinite and torsion free. M.~Brin and C.~Squier
showed that $F_2$ has no subgroups isomorphic to the free group of
rank 2 and satisfies no group laws~\cite{brin-s:nofree}. All normal
subgroups of $F_2$ contain the commutator, which is a simple
infinite group. The abelianization $F_2/[F_2,F_2]$ is ${\mathbb Z}
\times {\mathbb Z}$ (obvious from the presentation above). K.~Brown
and R.~Geoghegan singled out Thompson group $F_2$ as the first
example of a finitely presented torsion free group of type
$FP_\infty$ but not of type $FP$~\cite{brown-g:fp}. Thompson group
$F_2$ has a universal property with respect to homotopy
idempotents~\cite{freyd-h:hidempotents}. It is the group of order
preserving automorphisms of the free finitely generated algebra in
the variety of binary Cantor algebras (all finitely generated free
algebras in this variety are isomorphic, thus there is no notion of
a rank; see~\cite{higman:simple}). V.~Guba and M.~Sapir showed that
$F_2$ is the diagram group of the monoid presentation $\langle x
\mid x^2=x \rangle$~\cite{guba-s:diagram}. V.~Guba recently showed
that the Dehn function of $F_2$ is quadratic~\cite{guba:quadratic}
(this is exactly on the boundary between hyperbolic and
non-hyperbolic groups).

\begin{figure}[!ht]
\begin{center}
\includegraphics{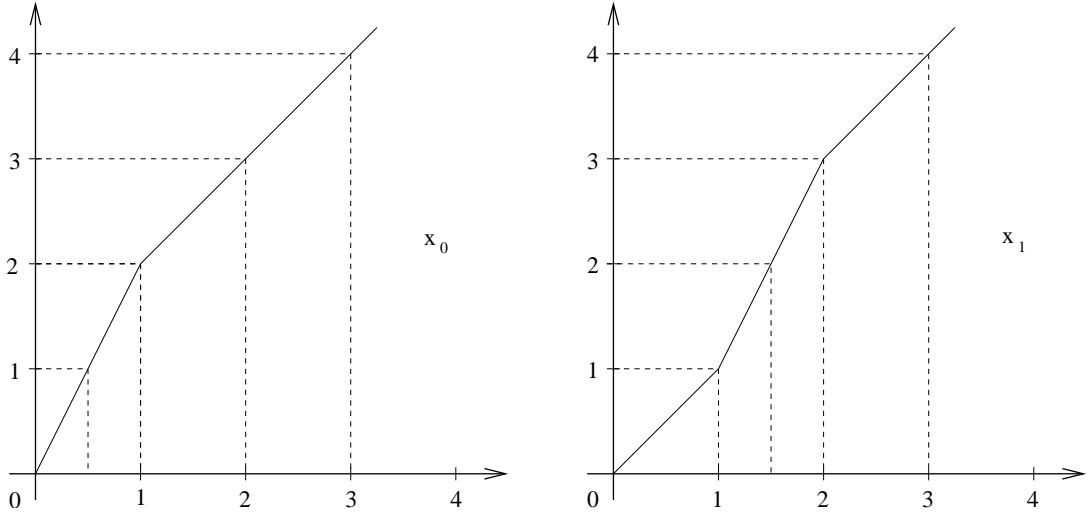}
\caption{The generators $x_0$ and $x_1$}\label{x0x1}
\end{center}
\end{figure}

On a very concrete level, the group $F_2$ may be realized as the
group of piecewise linear and order preserving homeomorphisms of the
unit interval $[0,1]$ such that all the slope breaks occur at dyadic
rational numbers and the slopes away from the finitely many breaks
are integer powers of 2 (this interpretation has been attributed to
Thurston). A closely related concrete realization is as the group
generated by the two piecewise linear homeomorphisms given in
Figure~\ref{x0x1} acting (on the left) on the interval $[0,\infty)$
.

The group $F_2$ was constructed by Thompson in 1965 in relation to
his study of questions in logic. The simplicity of the commutator
$[F_2,F_2]$ and the simplicity of two related finitely presented
groups, $T$ and $V$, were established by Thompson in his famous
unpublished notes~\cite{thompson:notes}. The groups $T$ and $V$ were
the first examples of finitely presented simple groups. A survey of
some properties of $F$ and the related simple groups $T$ and $V$ is
given in~\cite{cannon-f-p:thompson}.

It is known that $F_2$ is not elementary amenable, but it is not
known if it is amenable. J.~Belk and K.~Brown~\cite{belk-b:forest}
showed recently that the isoperimetric constant of $F_2$ is no
greater than $1/2$. The question of amenability of $F_2$ was raised
by R.~Geoghegan in~1979 and is one of the most interesting open
questions related to Thompson group $F_2$. The question of
amenability of $F_2$ can be related to the question of amenability
of the positive monoid $P_2$. It is shown by R.~Grigorchuk and
A.~Stepin in~\cite{grigorchuk-s:cancellation} that the positive
monoid $P_2$ is not left amenable (following the side convention we
use in this text) and that the group $F_2$ is amenable if and only
if the monoid $P_2$ is right amenable (the notion of amenability in
monoids requires attention to be paid to the side; left and right
amenability are the same in groups).

\section{Permutations, inversion sequences, linearized trees}

As seen from Section~\ref{fact-tamari} there is a long and fruitful
history of encoding permutations, trees (linearized or not) and
various integer sequences in terms of each other. We use this
section to establish a particular relation between trees,
permutations and integer sequences that is relevant to our
discussion and leads to a particular choice of a lattice congruence
${\sim}$ on $S_n$ defining the Tamari lattice $L_n=S_n/{\sim}$.

Consider a permutation $\sigma$ in $S_n$. Let
\[
 \mathsf{inv}_i(\sigma) = \# \{ j \mid 1 \leq j < \sigma^{-1}(i), \ \sigma(j) > i\ \}
\]
be the number of inversions of $\sigma$ that involve $i$ and a term
to the left of $i$ in $\sigma$. In other words,
$\mathsf{inv}_i(\sigma)$ counts the number of terms in $\sigma$ that
are larger than $i$ and are positioned to the left of the term $i$.
The sequence
\[ x(\sigma) = (\mathsf{inv}_1(\sigma), \mathsf{inv}_2(\sigma), \dots, \mathsf{inv}_n(\sigma) ) \]
is called the inversion sequence of $\sigma$. For, example, the
inversion sequence for the permutation $\sigma=32541$ is
$x(\sigma)=(4,1,0,1,0)$.

We will be thinking of sequences of natural numbers as elements of
the free monoid with basis ${\mathbb N}=\{0,1,\dots\}$. Let
$X=\{x_0,x_1,\dots\}$ and let $X^*$ be the free monoid on $X$. The
inversion sequence $(\mathsf{inv}_1(\sigma), \mathsf{inv}_2(\sigma),
\dots, \mathsf{inv}_n(\sigma))$ is then identified with the element
$x_{\mathsf{inv}_1(\sigma)}x_{\mathsf{inv}_2(\sigma)}\dots
x_{\mathsf{inv}_n(\sigma)}$ in $X^*$. Thus we have a map
\[ x: S_n \to X^* \]
that associates an $X$-word of length $n$ to any permutation in
$S_n$. Denote the image $x(S_n)$ by $X_n$. Then
\[
 X_n = \{ \ x_{i_1}x_{i_2} \dots x_{i_n} \mid 0 \leq i_j \leq n-j, \ j=1,\dots,n \ \}
\]
and $x:S_n \to X_n$ is bijective (see~\cite{stanley:ec1}). Let
\[ \pi : X_n \to S_n \]
be the inverse map of $x:S_n \to X_n$. One can try to write down
explicit formulae for $\pi(x)$, but it is more important for our
purposes to think of the following constructive way of calculating
the permutation $\pi(x)$ in $S_n$ corresponding to an inversion
sequence $x=x_{i_1}x_{i_2} \dots x_{i_n}$ in $X_n$. Start with $n$
empty slots $\underbrace{\_ \ \_ \ \dots \_ }_n$. Place 1 in such a
way that $i_1$ open slots are left to the left of it. Then place $2$
in such a way that $i_2$ open slots are left to the left of it.
Continue this procedure until a permutation in $S_n$ is obtained. In
other words, $\pi(x)$ is obtained in $n$ steps by placing, at step
$j$, the number $j$ in such a way that $i_j$ open slots are left to
the left of $j$, $j=1,\dots,n$. For example to calculate
$\pi(x_4x_1x_0x_1x_0)$ one starts with
\[ \begin{matrix} \_ & \_ & \_ & \_ & \_ & , \end{matrix} \]
then gets
\begin{equation}\label{steps}
\begin{matrix}
 \_ & \_  & \_ & \_ & 1 &, \\
 \_ & 2   & \_ & \_ & 1 &, \\
 3  & 2   & \_ & \_ & 1 &, \\
 3  & 2   & \_ & 4  & 1 &,
\end{matrix}
\end{equation}
and finally
\[ \begin{matrix} 3 & 2 & 5 & 4 & 1 &. \end{matrix} \]

Denote the set of linearized labeled rooted binary trees on $n$
interior vertices by $T_n$. Observe that the leafs $0,1,\dots,n$ are
usually drawn in line ordered from left to right by their labels.
Each par of consecutive leafs $(i-1,i)$, $i=1,\dots,n$, is called a
gap. Associate to each gap $(i-1,i)$ the last common vertex on the
unique paths from the root to leaf $i-1$ and leaf $i$. We say that
the associated interior vertex covers the gap (in order theoretic
terms this is just the join of the two leafs). In the example in
Figure~\ref{tree-lin} the correspondence between the gaps and the
interior vertices that cover them is given by
\[ \begin{matrix}
    3 & 2 & 5 & 4 & 1 \\
    \updownarrow & \updownarrow & \updownarrow & \updownarrow & \updownarrow \\
    (0,1) & (1,2) & (2,3) & (3,4) & (4,5)
    \end{matrix}.
\]
The correspondence between gaps and interior vertices that cover
them is bijective. The map
\[ \pi: T_n \to S_n \]
given by
\begin{equation}\label{pit}
 \pi(t)(i) = \text{ label of the interior vertex covering the gap } (i-1,i) \text{ in } t,
\end{equation}
for $i=1,\dots,n$ and a tree $t$ in $T_n$, is bijective. As already
observed, for our example from Figure~\ref{tree-lin} we have
$\pi(t)=32541$. We could express $\pi(t)$ without referring to the
gaps as follows. For $i=1,\dots,n$, $\pi(t)(i)$ is the linearization
label of the interior vertex visited at position $i$ in the
in-order. The reason we do not do this is that gaps will be relevant
in Section~\ref{partitions} when the discussion moves to trees of
higher degree.

A rooted binary tree with single interior vertex (and two leafs) is
called a caret. Each interior vertex determines a caret consisting
of the vertex itself and its two children. Each tree with $n$
interior vertices is composed of $n$ carets. Thus we can bijectively
associate gaps and carets in a labeled rooted binary tree. We also
say that the associated caret covers the corresponding gap.

\begin{figure}[!ht]
\begin{center}
\includegraphics{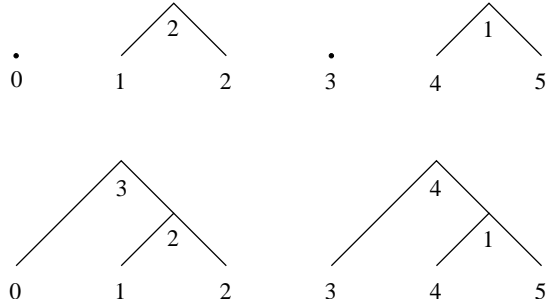}
\caption{Some intermediate steps in the construction of the
linearized tree $\tau(32541)$}\label{tree-lin2}
\end{center}
\end{figure}

We describe now the map
\[ \tau: S_n \to T_n, \]
which is inverse to $\pi: T_n \to S_n$, that associates a linearized
tree $\tau(\sigma)$ with $n$ interior vertices to a permutation
$\sigma$ in $S_n$. Start with $n+1$ leafs placed on a line and
labeled (from left to right) by $0,1,\dots,n$. In step $j$,
$j=1,\dots,n$, add a caret labeled by $j$ covering the gap
$(\sigma^{-1}(j)-1,\sigma^{-1}(j))$. For example, for the
permutation $\sigma=32541$ in the first step we add a caret labeled
by 1 covering the gap (4,5), after two steps we have two carets as
depicted in the top half of Figure~\ref{tree-lin2}, after 4 steps we
have 4 carets as in the bottom half of Figure~\ref{tree-lin2}
and in the last step we obtain the linearized tree in
Figure~\ref{tree-lin}.

Of course, the compositions $X_n \stackrel{\pi}{\to} S_n
\stackrel{\tau}{\to} T_n$ and its inverse $T_n \stackrel{\pi}{\to}
S_n \stackrel{x}{\to} X_n$ provide bijection between $X_n$ and $T_n$
for every $n$. It is actually possible to write down the
correspondence more directly without referring to $S_n$ as an
intermediate step, but a natural way to do this is to leave the
world of trees and inversion sequences and extend all considerations
to forests and arbitrary elements in the free monoid $X^*$. The
reason for this is already obvious when one notes that the
intermediate steps in Figure~\ref{tree-lin2} consists of forests
rather than trees. Similarly, the intermediate steps in
\eqref{steps} are not permutations in $S_n$. We want to develop a
language that will work with such intermediate steps. Additional
benefit of this extension is that, on the level of $X$-words, we
will work in the more natural environment of the full monoid $X^*$
rather than its submonoid $X_\infty = \cup_{n=0}^\infty X_n$.

Concatenation of inversion sequences of length $m$ and $n$ is an
inversion sequence of length $m+n$. Thus $X_\infty =
\cup_{n=0}^\infty X_n$ is indeed a submonoid of $X^*$. We can define
an operation on $S_\infty = \cup_{n=0}^\infty S_n$ that agrees with
the concatenation operation in $X_\infty$ and has the natural
extensions $x:S_\infty \to X_\infty$ and $\pi:X_\infty \to S_\infty$
as mutually inverse monoid isomorphisms. The operation is denoted by
$\div$ and defined as follows. If $\rho \in S_m$ and $\sigma \in
S_n$ then $\rho \div \sigma = \theta \in S_{m+n}$ is given by
\begin{equation}\label{operationS}
 \theta(i) =
    \begin{cases}
      \rho(i),        &1 \leq i \leq m \\
      \sigma(i-m)+ m, &m+1 \leq i \leq m+n
    \end{cases}.
\end{equation}
In other words, $\theta$ is produced by first increasing all the
terms of $\sigma$ by $m$ and then concatenating them to the right of
the terms of $\rho$. This operation on permutations will be called
interlacing (the same operation is used in~\cite{loday-r:hopftrees}
in the definition of product on the Hopf algebra $k[S_\infty]$).

\begin{figure}[!ht]
\begin{center}
\includegraphics{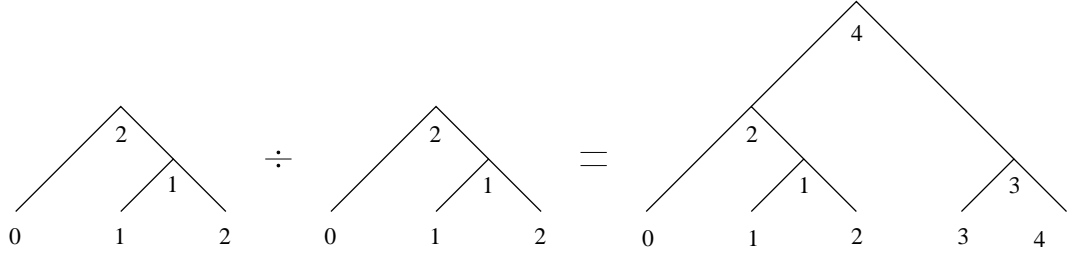}
\caption{Product of two linearized trees}\label{product-trees}
\end{center}
\end{figure}

Before we move on to forests and $X^*$ let us provide a definition
of the operation, also denoted by $\div$, on linearized trees in
$T_\infty = \cup_{n=0}^\infty T_n$ compatible with the concatenation
operation on $X_\infty$ and the interlacing operation $\div$ defined
in~\eqref{operationS} on $S_\infty$. The operation is performed by
stacking the second linearized tree on top of the first (hence the
notation $\div$). More precisely, for trees $r$ in $T_m$ and $s$ in
$T_n$ the product $t=rs$ is the tree $t$ in $T_{m+n}$ obtained by
deleting the leaf label 0 in $s$, increasing all other labels (both
in the interior and on the leafs) in $s$ by $m$, identifying the
leaf with deleted label in $s$ with the root of $r$ and declaring
the root od $s$ to be the root of $t$. An example is given in
Figure~\ref{product-trees}.

Thus right now we have three canonically isomorphic monoids, namely
$X_\infty$, $S_\infty$ and $T_\infty$, with the operations
concatenation, permutation interlacing and tree stacking,
respectively.

\section{Arbitrary sequences, linearized forests, $*$-permutations}

Since the free monoid $X^*$ does not need a special introduction, we
start by introducing the notion of linearized labeled binary rooted
forests. Such forests consists of countably many rooted binary
trees, only finitely many of which are non-trivial (the forest has
only finitely many carets). Furthermore, the roots are labeled in
bijective fashion by the numbers in ${\mathbb N}=\{0,1,2\dots\}$,
the leafs are also labeled in bijective fashion by the numbers in
${\mathbb N}$ and a leaf on one tree is labeled by a smaller number
than a leaf on another tree if and only if the same is true for
their corresponding roots. Note that vertices that are both a root
and a leaf have two labels, one as a root and one as a leaf, and
these labels may be different. Finally, if the number of carets is
$n$ then they are labeled bijectively by $1,2,\dots,n$ in such a way
that the labels on all paths from a root to a leaf are decreasing
(thus, again, the labeling of the interior vertices is compatible
with the order structure imposed by the forest structure). Denote
the set of linearized forests by $T^*$.

\begin{figure}[!ht]
\begin{center}
\includegraphics{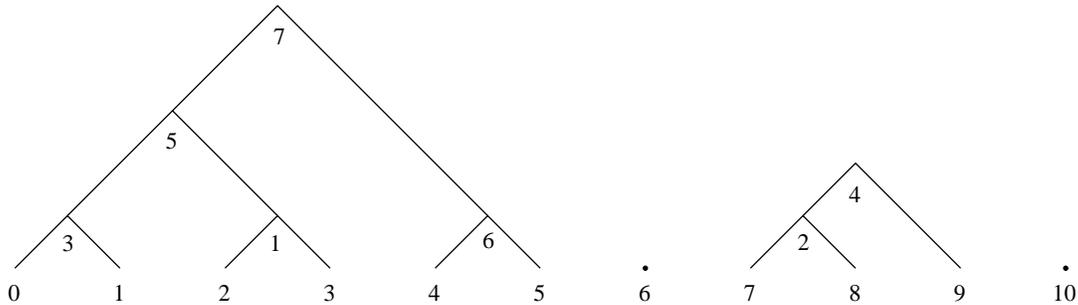}
\caption{A linearized forest}\label{forest-lin}
\end{center}
\end{figure}

We represent forests by diagrams of the type depicted in
Figure~\ref{forest-lin} in which it is assumed that the labeling of
both the roots and the leafs is done from left to right and all
trees that are not drawn are singletons labeled by higher numbers.
The labeling of the roots is usually left out, since it is
determined uniquely by the labeling of the leafs.
Note that leaf 6 is also labeled as root 1, while leaf 10 is also
labeled as root 3.

\begin{figure}[!ht]
\begin{center}
\includegraphics{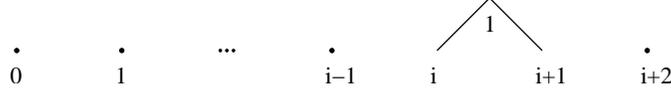}
\caption{The linearized forest corresponding to
$x_i$}\label{xi-lin}
\end{center}
\end{figure}

We define now a bijective map
\[ \tau: X^* \to T^*. \]
The forest corresponding to the word $x_{i_1}x_{i_2}\dots x_{i_n}$
in $X^*$ can be constructed inductively as follows. Start with the
trivial forest in which all trees are singletons. Throughout the
whole construction the leafs and their labels are left unchanged.
All that happens in the process is that we add labeled carets and
relabel the roots. For a letter $x_i$ the corresponding linearized
forest is given in Figure~\ref{xi-lin}.
Note that the root labels to the right of leaf $i+1$ do not agree
any longer with the leaf labels (the root labels are  smaller by 1).
If $f_{n-1}$ is the forest representing $x_{i_1}x_{i_2}\dots
x_{i_{n-1}}$ construct the forest $f_n$ representing
$x_{i_1}x_{i_2}\dots x_{i_n}$ by adding a caret, labeled by $n$,
connecting root $i_n$ and root $i_n+1$. The newly added root gets
root label $i_n$, old roots $i_n+2$ and higher get their labels
decreased by 1 and the old roots $i_n$ and $i_n+1$ loose their root
labels (they are not roots any longer). For example, the forest in
Figure~\ref{forest-lin} corresponds to the word
$x_2x_6x_0x_5x_0x_1x_0$.

Conversely, the inverse map
\[ x: T^* \to X^* \]
can be understood as follows. The linearization part of the labeling
of the forest $f$ gives a recipe for constructing $f$ starting from
the trivial forest. Namely, first add the caret labeled by 1, then
the one labeled by 2, and so on until the caret labeled by $n$ is
added. In each step $j$ we record the label $i_j$ of the root that
the newly introduced caret labeled by $j$ uses as the left leaf. The
corresponding element of $X^*$ is then $x_{i_1}x_{i_2}\dots
x_{i_n}$. One can check that this procedure applied to the tree in
Figure~\ref{forest-lin} does indeed yield the word
$x_2x_6x_0x_5x_0x_1x_0$.

\begin{figure}[!ht]
\begin{center}
\includegraphics{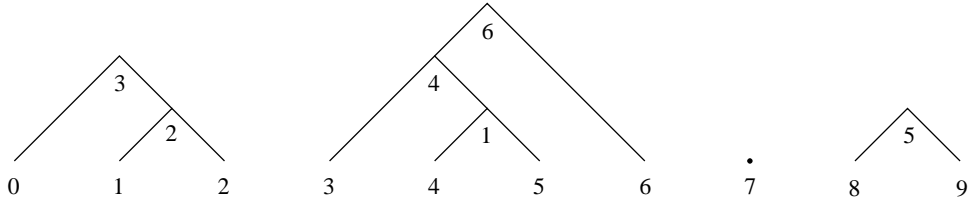}
\caption{Product of two linearized forests}\label{product-forest}
\end{center}
\end{figure}

Forests can be multiplied in a way compatible with the concatenation
multiplication in $X^*$ as follows. The operation, still denoted by
$\div$, is performed by stacking the second forest on top of the
first. Namely, the product of the forests $f$ and $g$ is the forest
$h$ obtained by, first increasing all labels on the carets of $g$ by
$m$, where $m$ is the number of carets in $f$, then identifying root
$i$ in forest $f$ with leaf $i$ in forest $g$ and then deleting
their root/leaf labels, correspondingly. The leafs of $h=fg$ are the
leafs of $f$ and the roots of $h=fg$ are the roots of $g$. For
example, the product of the forest in the bottom half of
Figure~\ref{tree-lin2} and the forest in the top half of the same
figure is the linearized forest in Figure~\ref{product-forest}.

We now turn to the world of permutations. The intermediate steps in
\eqref{steps} provide motivation for the following definition. A
$*$-permutation of length $n$ is a map $\sigma: {\mathbb N}^+ \to
\{1,2,\dots,n\}\cup \{*\}$, where ${\mathbb N}^+=\{1,2,\dots\}$,
such that the inverse image of each element in $\{1,2,\dots,n\}$ is
a singleton. In other words, $\sigma$ is an infinite sequence such
that each of $1,\dots,n$ appears exactly once as a term, and the
rest of the terms are $*$'s. Denote the set of $*$-permutations by
$S^*$.

A bijective correspondence
\[ \pi: T^* \to S^* \]
between linearized forests and $*$-permutations is defined as
follows. Each gap $(i-1,i)$, $i=1,2,\dots$ in a linearized forest
$f$ is either covered by a caret or is not covered by a caret (the
latter happens in case the two vertices defining the gap come from
different trees). Define
\[ \pi(f)(i) =
    \begin{cases}
       \text{label of the caret covering the gap }(i-1,i)\text{ in }f, & \text{if a cover exists}\\
       * & \text{otherwise}
    \end{cases}.
\]
For example, the $*$-permutation corresponding to the forests in
Figure~\ref{product-forest} and Figure~\ref{forest-lin} are
$32*416**5$ and $35176**24*$, respectively. Note that we agree
(sometimes) not to write (some of) the stars to the right of the
last non-$*$ symbol. The non-$*$ symbols will be called concrete
symbols in the rest of the text.

The inverse map
\[ \tau:S^* \to T^* \]
is, just as its restriction to $S_\infty$, simply defined by drawing
the forest caret by caret. For a $*$-permutation of length $n$, in
step $j$, $j=1,\dots,n$, add a caret labeled by $j$ covering the gap
$(\sigma^{-1}(j)-1,\sigma^{-1}(j))$ and relabel the roots
accordingly. At the end of the procedure only $n$ carets will be
drawn and the rest of the gaps are uncovered and correspond to
$*$'s.

The mutually inverse maps
\[ \pi:X^* \to S^* \quad\text{and}\quad x:S^* \to X^* \]
are defined in the same fashion as before. If $\sigma$ has $n$
concrete symbols the corresponding word $x(\sigma)$ in $X^*$ is
$x_{\mathsf{inv}_1(\sigma)}x_{\mathsf{inv}_2(\sigma)}\dots
x_{\mathsf{inv}_n(\sigma)}$, with added stipulation that any
occurrence of a $*$ to the left of $j$ in $\sigma$ is counted in
$\mathsf{inv}_j(\sigma)$. In other words, for all our purposes, $*$
is considered to be larger than any concrete symbol (symbol from
${\mathbb N}^+$). In the other direction, given a word $x$ of length
$n$ over $X$, once the symbols $1,2,\dots,n$ are placed by leaving
an appropriate number of open slots as prescribed by $x$, the rest
of the open slots are filled with $*$'s. For example, $32*416**5$
and $35176**24$ correspond to the $X$-words $x_4x_1x_0x_1x_4x_1$ and
$x_2x_6x_0x_5x_0x_1x_0$, respectively.

To complete the picture we define the operation, still denoted by
$\div$, on $S^*$ that agrees with concatenation on $X^*$ and
stacking of linearized trees in $T^*$. The operation is performed by
interlacing $*$-permutations. One can write down definite formulae,
but the operation is more easily understood as a process. For
$*$-permutation $\rho$ with $m$ concrete symbols and $*$-permutation
$\sigma$ with $n$ concrete symbols the product $\theta = \rho
\sigma$ is obtained as follows. First increase the concrete symbols
in $\sigma$ by $m$. Then interlace $\sigma$ into $\rho$ by placing
the $i$-th symbol of $\sigma$ in the position of the $i$-th star in
$\rho$. For example $(25*31**4) (3*1*2)= 25831*64*7$. The following
diagram may be helpful for imagining the process
\begin{center}
\begin{tabular}{cccccccccccccccccccl}
 ( & 2 & 5 & * & 3 & 1 & * & * & 4 & ) &$\div$& ( & 3 & * & 1 & * & 2 & )  &  = \vspace{3mm}\\
   &   &   &   &   &   &   &   &   &   &      &   & 8 & * & 6 & * & 7 &    &    &lift in value\\
   & 2 & 5 & * & 3 & 1 & * & * & 4 &   &      &   &   &   &   &   &   &    &  = &and literally \vspace{3mm}\\
   &   &   & 8 &   &   & * & 6 &   &  * & 7   &   &   &   &   &   &   &    &    &place above stars \\
   &   &   & $\downarrow$ &   &   & $\downarrow$  & $\downarrow$ & &    $\downarrow$ & $\downarrow$ \\
   & 2 & 5 & * & 3 & 1 & * & * & 4 & * & *    &   &   &   &   &   &   &     & = \vspace{3mm} \\
   & 2 & 5 & 8 & 3 & 1 & * & 6 & 4 & * & 7    &   &   &   &   &   &   &     &    &lower in place of stars
\end{tabular}
\end{center}

\section{Tamari congruence induced by
de-linearization}\label{delinearization}

At this moment we have three monoids $X^*$, $T^*$ and $S^*$ related
by canonical isomorphisms. Since $X^*$ is free so are $T^*$ and
$S^*$ and we may be disappointed that all that happened so far is
that we obtained two strange copies of the free monoid $X^*$ of
countable rank -- namely $T^*$ with a free basis consisting of
linearized forests $t_i$ as in Figure~\ref{xi-lin}, $i\in {\mathbb
N}$, and $S^*$ with a free basis consisting of
$s_i=\underbrace{**\dots
*}_i 1$, $i \in {\mathbb N}$.

There are at least two ways to motivate what comes next.

One is to observe that the multiplication rule $\div$ on $T^*$ does
not essentially depend on the linearization part of the labeling of
the involved forests. This labeling is just carried around and
adjusted here and there by increasing labels accordingly, but
nothing in the definition depends on it. This means that the
equivalence relation ${\sim}$ on $T^*$ obtained by dropping the
labels on interior vertices is not only equivalence on $T^*$ but it
is also a monoid congruence.

\begin{proposition}
The equivalence ${\sim}$ is a congruence on the monoid $(T^*,\div)$.
\end{proposition}

Another way to motivate the introduction of ${\sim}$ is as an
extension of a well known connection between permutations and their
linearized trees obtained when the linearization is striped away. In
that case, several permutations correspond to the same labeled
rooted binary tree. It is known that there are Catalan number $C_n=
\frac{1}{n+1}\binom{2n}{n}$ labeled rooted binary trees on $n$
interior vertices. Thus the $n!$ permutations in $S_n$ are split
into $C_n$ classes of equivalent permutations. We fix the
equivalence classes obtained in this process as classes defining the
Tamari congruence on $S_n$.

Formally, for any forest $t$ in $T^*$ define $\widetilde{t}$ to be
the forest obtained when the labeling on all interior vertices is
deleted. Define an equivalence on $T^*$ by
\[ r \sim t \Leftrightarrow \widetilde{r} = \widetilde{t} \]
and, by use of the corresponding bijections, define the induced
equivalences on $S^*$ by
\[ \rho \sim \sigma \Leftrightarrow \tau(\rho) \sim \tau(\sigma) \Leftrightarrow \widetilde{\tau(\rho)} = \widetilde{\tau(\sigma)} \]
and on $X^*$ by
\[ u \sim v \Leftrightarrow \tau(u) \sim \tau(v) \Leftrightarrow \widetilde{\tau(u)} = \widetilde{\tau(v)}. \]
Thus we have a monoid congruence ${\sim}$ on $X^*$, $S^*$ and $T^*$
and we want to understand the corresponding factor monoid.

\begin{figure}[!ht]
\begin{center}
\includegraphics{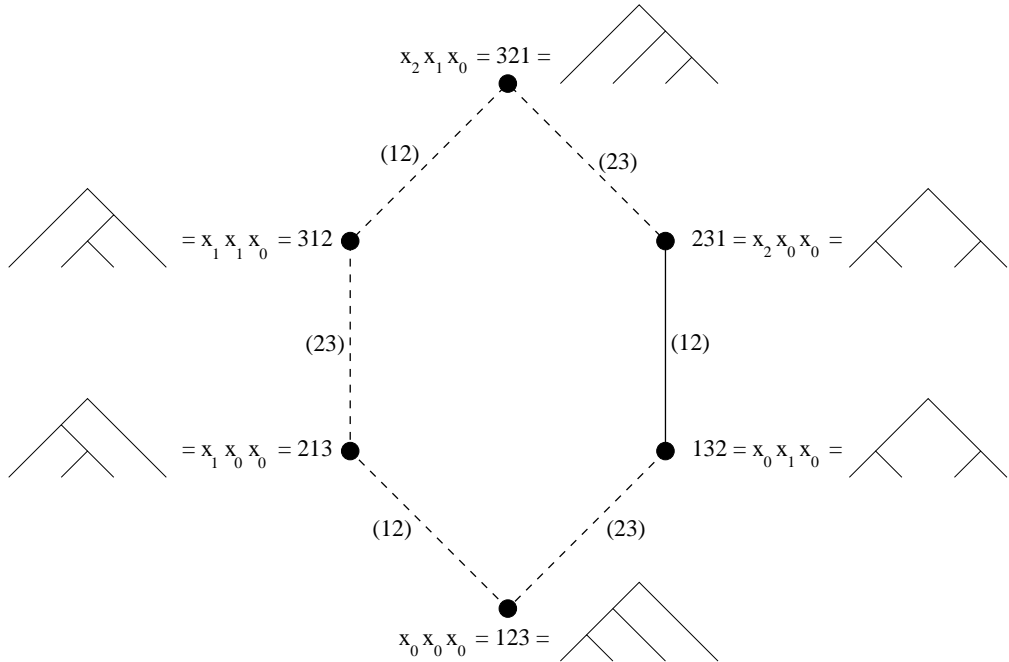}
\caption{$S_3$ and its associated words and trees}\label{s3all}
\end{center}
\end{figure}

The diagram in Figure~\ref{s3all} depicts the situation for $S_3$.
The hexagon in the middle is the (left) Cayley graph of $S_3$ as
Coxeter group of type $A_2$ generated by the standard generating
set $\{(12),(23)\}$. It is drawn in a way that represents the
Hasse diagram of the (left) weak Bruhat order on $S_3$.
The only edge drawn as a full line in the right half of the diagram
indicates that the corresponding permutations (or words) are to be
identified in $S_3/{\sim}$ (or in $X_3/{\sim}$), since they
correspond to the same tree. Therefore $132 \sim 231$, $x_0x_1x_0
\sim x_2x_0x_0$ and $S_3/{\sim}$ and $X_3/{\sim}$ have 5 elements
(which is the Catalan number $C_3$).

Observe that the restriction of ${\sim}$ to $S_\infty$ is a
congruence on the monoid $S_\infty$ (again, this is because when the
labeling is stripped in $T_\infty$ the operation $\div$ is not
affected).

\begin{proposition}\label{free}
The monoid $S_\infty/{\sim} = T_\infty/{\sim} = X_\infty/{\sim}$ is
free.
\end{proposition}
\begin{proof}
Indeed, any non-trivial tree for which the left subtree is trivial
is indecomposable in $T_\infty/{\sim}$ (it cannot be written as a
product of two or more nontrivial trees). Moreover, every tree in
$T_\infty/{\sim}$ has a unique decomposition as a product of such
indecomposable trees. Thus $T_\infty/{\sim}$ is free.
\end{proof}

This is perhaps a reason why the connection to the Thompson monoid
is not immediately obvious. For a researcher in combinatorics there
seem to be no particular gain in constructing free monoids using
strange definitions of products of permutations. On the other hand,
looking from Thompson monoid $P_2$ point of view, the connection to
Tamari lattices on $S_n$ is not immediately obvious since on its
basic level, working only with ordinary permutations before the
expansion to $S^*$, the information is encoded in a certain (not
particularly distinguished) free submonoid of $P_2$ that seemingly
does not demand any particular attention (there are plenty of free
monoids inside $P_2$).

However, we want to understand $T^*/{\sim}=X^*/{\sim} = S^*/{\sim}$
and this is where the interesting things happen. Our understanding
of the weak Bruhat order and Tamari lattice congruences on $S_n$ can
still be used in this extended situation.

Note that every element in $S^*$ has finite congruence class. This
is because there are only finitely many ways to linearize a forest
with finitely many carets.

Let $S_n'$ be the set of $n!$ different $*$-permutations with $n$
concrete (non-$*$) terms and $*$'s in some fixed positions. More
precisely, let us say that there are exactly $k$ blocks of
consecutive positions in which concrete symbols appear (any two
blocks are separated by some finite number of $*$'s). Let the sizes
of the concrete blocks, from left to right, be $m_1,\dots,m_k$ and
let the last concrete symbol appears at position $z$. Clearly, all
the elements related by ${\sim}$ to a $*$-permutation in $S_n'$ are
also in $S_n'$. We will describe the congruence classes on $S_n'$ in
terms of the congruence classes on $S_{m_1},\dots,S_{m_k}$.

There is a canonical correspondence $S_n' \leftrightarrow S_n$
obtained by removal/insertion of $*$'s in appropriate positions.
This enables us to induce the weak Bruhat order on $S_n'$. We write
$(i \ i+1) \circ \sigma$ for the $*$-permutation obtained from
$\sigma$ when $i$ and $i+1$ exchange their positions in $\sigma$.
Extending to $S_n$ this defines a (left) action of $S_n$ on $S_n'$.
Technically speaking, the Hasse diagram of the weak Bruhat order
induced on $S_n'$ is not the Cayley graph of $S_n$ but rather the
Schreier graph of the action of $S_n$ on $S_n'$ with respect to the
standard generating set $\{(12),\dots,(n-1 \ n)\}$, but these two
graphs are canonically isomorphic and we borrow the terminology from
$S_n$ and use it on $S_n'$. In particular, we keep the notation
$\preceq$ for the weak Bruhat order in the extended sense.

If $\rho=(i \ i+1) \circ \sigma$ in $S_n'$ is obtained from
$\sigma$ by interchanging $i$ and $i+1$ it is still valid to say
that the lengths of $\rho$ and $\sigma$ differ by 1. Moreover,
$\rho$ covers $\sigma$ in the weak Bruhat order if and only if $i$
is to the left of $i+1$ in $\sigma$ and it is directly below
$\sigma$ in the other case (this says that the weak Bruhat order
is compatible with the lexicographic order on $S_n'$).

\begin{figure}[!ht]
\begin{center}
\includegraphics{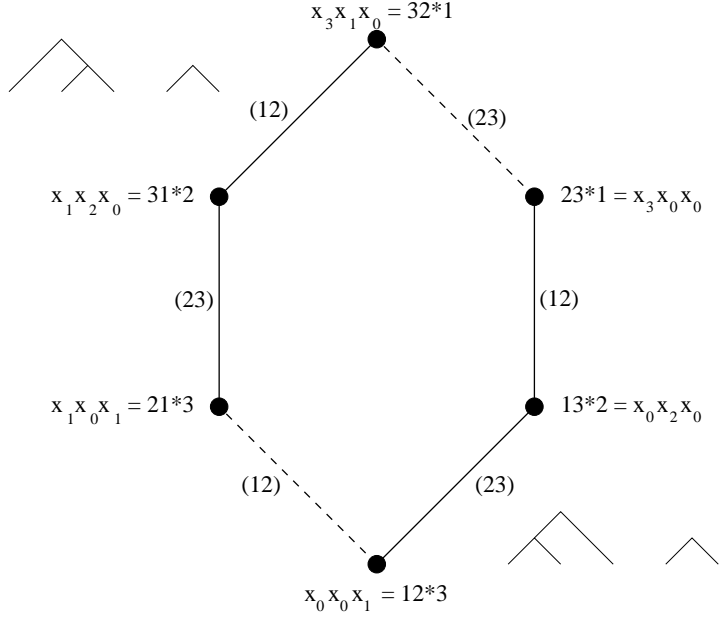}
\caption{Equivalence classes in $S_3'$ corresponding to the block
pattern $\_ \ \_ * \_$}\label{s3prime}
\end{center}
\end{figure}

We note that there is a very important difference. Namely the
${\sim}$ classes on $S_n$ and on $S_n'$ are not the same. For
example, $S_3$ has five classes, while $S_3'$ corresponding to the
block pattern $\_ \ \_ * \_$ has only two. The two classes in $S_3'$
are indicated in Figure~\ref{s3prime} as the components connected by
edges drawn as full lines.

The first thing we want to know is under what conditions two
neighbors in the weak Bruhat order correspond to the same forest.
We recall the explicit agreement that $*$ is larger than any
concrete symbol.

\begin{proposition}\label{move-up}
Let $\rho = ( i \ i+1) \circ \sigma$ in $S_n'$. Then the forests
$\tau(\rho)$ and $\tau(\sigma)$ are equal if and only if there
exists a term between $i$ and $i+1$ in $\sigma$ (and therefore in
$\rho$ as well) that is larger than $i+1$.
\end{proposition}
\begin{proof}
Exchanging the labels $i$ and $i+1$ in a linearized forest $f$ leads
to another linearized forest if and only if the interior vertices
$i$ and $i+1$ are not comparable with respect to the order induced
by the forest $f$. This happens exactly when $i$ and $i+1$ belong to
different trees in the forest, in which case there is a $*$ between
$i$ and $i+1$ in $\pi(f)$, or $i$ and $i+1$ belong to the same tree,
but are both descendants of a vertex $j$, in which case $j>i+1$ and
$j$ is between $i$ and $i+1$ in $\pi(f)$.
\end{proof}

This means that for every $\sigma$ in $S_n'$ we may go up step by
step in the weak Bruhat order until we reach a $*$-permutation
that has no occurrence of $\dots i \dots j \dots i+1 \dots$, with
$j>i+1$. This leads us to pattern avoiding considerations.

We say that $\overline{132}$ occurs in $\sigma$ in $S_n'$ if there
are three indices $1 \leq a < b < c \leq z$ such that $\sigma(a) + 1
= \sigma(c) < \sigma(b)$ (recall that $z$ is the index of the last
concrete symbol in the $*$-permutations in $S_n'$). This allows the
possibility that the middle symbol $\sigma(b)$ is a $*$, but the
other two symbols must be concrete. If $\overline{132}$ does not
occur in $\sigma$ the $*$-permutation is $\overline{132}$-avoiding.
Similarly, we say that 132 occurs in a $*$-permutation $\sigma$ in
$S_n'$ (or in any sequence over a linearly ordered set) if there are
three indices $1 \leq a < b < c \leq z$ such that $\sigma(a) <
\sigma(c) < \sigma(b)$. If 132 does not occur in $\sigma$ the
permutation is 132-avoiding. Once again, the definition implies that
the middle symbol $\sigma(b)$ may be a $*$ but the other two symbols
involved must be concrete.

It is easy to observe that $\sigma$ is a 132-avoiding
$*$-permutation if and only if all concrete terms to the left of
any occurrence of $*$ in $\sigma$ are larger than all concrete
terms to the right of the same occurrence of $*$ and each concrete
block in $\sigma$ satisfies the 132-avoiding constraint by itself.

\begin{proposition}\label{132=132}
A $*$-permutation $\sigma$ in $S_n'$ is $132$-avoiding if and only
if it is $\overline{132}$-avoiding.
\end{proposition}
\begin{proof}
Any occurrence of $\overline{132}$ in $\sigma$ is also an
occurrence of $132$.

For the converse, assume $132$ appears in $\sigma$. Let $\dots i
\dots k \dots j \dots$ be an occurrence of $132$ with minimal
difference $j-i$. If $i+1$ is to the right of $k$ we have an
appearance of $\overline{132}$. But $i+1$ cannot appear to the
left of $k$ because of the minimality in the choice of $i$ and
$j$.
\end{proof}

The above observations lead to the following proposition.

\begin{proposition}\label{unique-up}
For every $*$-permutation $\sigma$ in $S_n'$ there exists a unique
$132$-avoiding $*$-permutation $\overline{\sigma}$ in $S_n'$ such
that $\sigma \sim \overline{\sigma}$. Moreover $\sigma \preceq
\overline{\sigma}$ and $\rho \sim \sigma$ if and only if
$\overline{\rho} = \overline{\sigma}$.
\end{proposition}
\begin{proof}
The existence is clear (just go up step by step as long as
possible without changing the forest).

Everything else follows from counting arguments.

There are exactly $C_{m_1}C_{m_1}\dots C_{m_k}$ ordered $k$-tuples
of labeled rooted binary trees with $m_1,\dots,m_k$ carets,
respectively.

As for $132$-avoiding $*$-permutations, the symbols in each
concrete block are uniquely determined by the condition that all
symbols in a block to the left of some other block must be larger
than the symbols in the other block. Since there are exactly
$C_{m_i}$ $132$-avoiding arrangements of $m_i$ symbols (see the
Catalan addendum in~\cite{stanley:ec2}) we get that there are
$C_{m_1}C_{m_1}\dots C_{m_k}$ $132$-avoiding $*$-permutations in
$S_n'$. Thus every equivalence class in $S_n'$ must have exactly
one $132$-avoiding representative.
\end{proof}

In a completely analogous way, starting from any $*$-permutation
$\sigma$ in $S_n'$ we may move step by step down in the weak
Bruhat order by preserving the corresponding forest as long as we
see an occurrence of $\dots i+1 \dots j \dots i \dots$, with $j
>i+1$. Define $231$-avoiding $*$-permutations as $*$-permutations in which
there are no three indices $1 \leq a < b < c \leq z$ such that
$\sigma(c) < \sigma(a) < \sigma(b)$.

\begin{proposition}\label{unique-down}
For every $*$-permutation $\sigma$ in $S_n'$ there exists a unique
$231$-avoiding $*$-permutation $\underline{\sigma}$ in $S_n'$ such
that $\sigma \sim \underline{\sigma}$. Moreover $\underline{\sigma}
\preceq \sigma$ and $\rho \sim \sigma$ if and only if
$\underline{\rho} = \underline{\sigma}$.
\end{proposition}

It is known that $L_n=S_n/{\sim}$ is a lattice, known as Tamari
lattice. In other words, it is known that ${\sim}$ is a lattice
congruence on $S_n$. We claim that $S_n'/{\sim}$ is also a lattice,
i.e. ${\sim}$ is a lattice congruence on the weak Bruhat order
lattice on $S_n'$. This essentially follows from the fact that
certain maps $S_n \to S_m$, where $S_m$ is parabolic subgroup of
$S_n$ (as Coxeter groups) are lattice homomorphisms on the
corresponding weak Bruhat order lattices.

Define the flattening of a sequence $e_1,e_2,\dots,e_m$ of
distinct integers (or members of any linear order) as the unique
permutation $\sigma$ in $S_m$ such that $e_i < e_j$ if and only if
$\sigma(i) < \sigma(j)$.

\begin{lemma}\label{left-right}
The map $S_n \to S_m$ that maps a permutation $\sigma$ to the
flattening of $\sigma(p+1)\sigma(p+2)\dots\sigma(p+m)$ is a
surjective lattice homomorphism.
\end{lemma}
\begin{proof}
The statement is equivalent to the following. Consider $S_n$ and
$S_m$ under their right weak Bruhat order. Define a map $\alpha: S_n
\to S_m$ as follows. For $\sigma$ in $S_n$, let $\alpha(\sigma)$ be
the permutation in $S_m$ obtained when all terms in $\sigma$ except
for the terms $p+1,\dots,p+m$ are deleted and then flattened
(decreased by $p$). The equivalence comes by applying the inverse to
the elements in $S_n$. In the left weak Bruhat order we are
concerned with the positions $p+1,\dots,p+m$, while after inversion
takes place and we land in the right weak Bruhat order we are
concerned with the terms $p+1,\dots,p+m$. With respect to the right
Bruhat order, the map $\alpha$ is the surjective lattice
homomorphism $S_n \to S_m$ corresponding to the parabolic subgroup
$S_m$ generated by the $m-1$ reflections $\{(p+1 \ p+2), \dots,
(p+m-1 \ p+m)\}$ (see~\cite{reading:congruences} for example).
\end{proof}

\begin{theorem}
The equivalence ${\sim}$ is a lattice congruence on $S_n'$. Moreover
\[ S_n'/{\sim} \cong S_{m_1}/{\sim} \times \dots \times S_{m_k}/{\sim} \]
as lattices.
\end{theorem}
\begin{proof}
Patching together $k$ parabolic homomorphisms as in
Lemma~\ref{left-right} we get a surjective lattice homomorphism
$S_n' \to S_{m_1} \times \dots \times S_{m_k}$, which can then be
composed further to get a lattice homomorphism $S_n' \to
S_{m_1}/{\sim} \times \dots \times S_{m_k}/{\sim}$. We claim that
${\sim}$ is the kernel of this homomorphism. Recall that the
concrete terms of a $*$-permutation $\sigma$ in $S_n'$ just indicate
in what order the carets are added in the forest $\tau(\sigma)$,
which consists of $k$ trees with $m_1,\dots,m_k$ carets,
respectively. In particular, it is clear that the $i$-th
de-linearized tree that corresponds to the $i$-th concrete block
depends only on the ${\sim}$ class of the flattening of the
corresponding block (the gaps in the numbers before the flattening
correspond to carets added in the other trees of the forest). Thus
two $*$-permutations correspond to the same forest if and only if
the corresponding flattenings in each block are ${\sim}$ related and
the relation ${\sim}$ on $S_n'$ is indeed the kernel of the
surjective homomorphism $S_n' \to S_{m_1}/{\sim} \times \dots \times
S_{m_k}/{\sim}$.
\end{proof}

\begin{corollary}\label{intervals}
Each ${\sim}$ conjugacy class of $*$-permutations in $S_n'$ is an
interval in the weak Bruhat order and is a union of several Tamari
congruence classes of $S_n$ (after identification of $S_n'$ and
$S_n$).

For each class, the top bound of the interval is a $132$-avoiding
$*$-permutation and the bottom bound is a $231$-avoiding
$*$-permutation.
\end{corollary}
\begin{proof}
Congruence classes in finite lattices are always intervals.

Proposition~\ref{unique-up} and Proposition~\ref{unique-down} show
that the top and the bottom must be $132$-avoiding and
$231$-avoiding $*$-permutations, respectively.

Finally, if $\rho= (i \ i+1) \circ \sigma$ and $\rho \sim \sigma$ in
$S_n$, then there exists a term between $i$ and $i+1$ larger than
$i+1$. After $*$'s are placed in appropriate places to land in
$S_n'$ it is still correct that there is a larger term between $i$
and $i+1$. Thus the corresponding $*$-permutations are also related.
This shows that each Tamari class of $S_n$ is included in a ${\sim}$
class of $S_n'$.
\end{proof}

Thus we have a thorough understanding of the equivalence classes
in $S^*$. We translate now this understanding to $X^*$.

We know that we can connect any two equivalent $*$-permutations by
several steps involving transpositions, with the extra constraint
that when we apply $(i \ i+1)$ some term between $i$ and $i+1$ must
be larger than $i+1$. Here is the corresponding statement in the
$X^*$ world.

\begin{proposition}\label{move-upX}
Let $\rho = (j \ j+1) \circ \sigma $, $\sigma \preceq \rho$ and
$\rho \sim \sigma$ in $S^*$. Further let $x(\sigma) = \dots
x_{\mathsf{inv}_j(\sigma)}x_{\mathsf{inv}_{j+1}(\sigma)} \dots =
\dots x_m x_n \dots$~. Then $m<n$ and $x(\rho)$ can be obtained from
$x(\sigma)$ by applying the substitution
\[ x_mx_n \to x_{n+1}x_m \]
at positions $j$ and $j+1$ in $x(\sigma)$.
\end{proposition}
\begin{proof}
By Proposition~\ref{move-up}, $\sigma = \dots j \dots k \dots j+1
\dots$, where $k$ is a term larger than $j+1$ (possibly a $*$).
This immediately means that there are more inversions to the left
of $j+1$ than to the left of $j$, i.e., $m<n$. The exchange of $j$
and $j+1$ causes move upwards in the weak Bruhat order
\begin{center}
\begin{tabular}{CCCCCCCCC}
  \rho   & = & \dots & j+1 & \dots & k & \dots & j  & \dots \vspace{3mm}\\
         &   &       &     &       &  \uparrow \ (j \ j+1) \vspace{3mm} \\
  \sigma & = & \dots & j   & \dots & k & \dots & j+1 & \dots
\end{tabular}
\end{center}

After $j$ and $j+1$ exchange their positions, no inversion numbers
other than those at position $j$ and $j+1$ in $x(\sigma)$ can
possibly be affected. We have $\mathsf{inv}_{j+1}(\rho) =
\mathsf{inv}_j(\sigma) = m$ and $\mathsf{inv}_j(\rho) =
\mathsf{inv}_{j+1}(\sigma)+1 = n+1$. The extra 1 in
$\mathsf{inv}_j(\rho)$ comes from the fact that now $j+1$ is to the
left of $j$ and should be counted as extra inversion. Thus
\begin{center}
\begin{tabular}{CCCL}
  x(\rho)   & = & \dots \ x_{n+1} x_m \ \dots \vspace{3mm}\\
            &   &   \uparrow & \text{ at positions }(j,j+1) \vspace{3mm} \\
  x(\sigma) & = & \dots \ x_m x_n  \ \dots
\end{tabular}
\end{center}
\end{proof}

A converse to the previous proposition holds.

\begin{proposition}
Let $x = \dots x_{i_j}x_{i_{j+1}} \dots = \dots x_m x_n ...$, with
$m<n$, and let $y$ be obtained from $x$ by applying the substitution
$x_mx_n \to x_{n+1}x_m$. Then $\pi(y) = (j \ j+1) \circ \pi(x)$,
$\pi(x) \preceq \pi(y)$ and $\pi(x) \sim \pi(y)$.
\end{proposition}

\begin{corollary}\label{xsimy}
Let $x = \dots x_{i_j}x_{i_{j+1}} \dots = \dots x_m x_n ...$, with
$m<n$, and let $y$ be obtained from $x$ by applying the substitution
$x_mx_n \to x_{n+1}x_m$. Then $x \sim y$.
\end{corollary}

We can now prove that the monoid $X^*/{\sim}$ is Thompson's monoid
$P_2$.

\begin{theorem}
The congruence ${\sim}$ on $X^*$ is generated by
\[ x_ix_j \sim x_{j+1}x_i, \]
for all pars of non-negative integers $i$ and $j$ with $i <j$.

In other words, $X^*/{\sim}$ is equal to Thompson's monoid $P_2$,
given by the presentation
\[
 Mon \langle \ x_0,x_1,x_2,\dots \mid x_i x_j = x_{j+1}x_i, \text{ for }i<j \ \rangle.
\]
\end{theorem}
\begin{proof}
It is clear that $x_ix_j \sim x_{j+1}x_i$ does hold, for $i <j$, in
$X^*$ (this is a special case of Corollary~\ref{xsimy}). Another way
to see this is to realize that the two words $x_ix_j$ and
$x_{j+1}x_i$ correspond to the two ways to linearize (i.e.\ to draw)
the forest in Figure~\ref{xixj}
\begin{figure}[!ht]
\begin{center}
\includegraphics{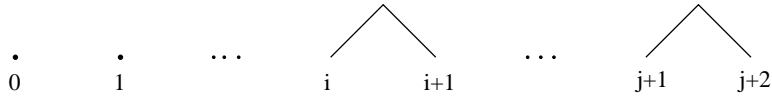}
\caption{$x_ix_j \sim x_{j+1}x_i$}\label{xixj}
\end{center}
\end{figure}

On the other hand, the relations $x_ix_i \sim x_{j+1}x_i$, for
$i<j$, are sufficient. This is because any two equivalent
$*$-permutations can be related by a sequence of applications of
appropriate transpositions (say by going up in the weak Bruhat order
on $S_n'$ and reaching the same $132$ avoiding $*$-permutation). By
Proposition~\ref{move-upX} this translates to a sequence of
applications of relations of the type $x_ix_j \sim x_{j+1}x_i$, for
$i<j$.
\end{proof}

\begin{corollary}
Every class of equivalent elements in $P_2$ corresponds to an
interval in the weak Bruhat order on $*$-permutations. Moreover, the
top always corresponds to a $132$-avoiding and the bottom to a
$231$-avoiding $*$-permutation.
\end{corollary}

The following proposition characterizes the words over $X$ that
correspond to the endpoints of weak Bruhat order intervals that
determine the classes of equivalent elements in Thompson's monoid
$P_2$.

\begin{proposition}
A $*$-permutation is $132$-avoiding if and only if the
corresponding word $x(\sigma)$ has non-increasing indices.

A $*$-permutation is $231$-avoiding if and only if the
corresponding word $x(\sigma)$ has no decrease of an index larger
than 1.
\end{proposition}

As a corollary we obtain the two well known normal forms on $P_2$.

\begin{corollary}
Every class of equivalent elements in $P_2$ has two normal forms.
One is a word with non-increasing indices and the other is a word
whose indices never decrease by more than 1 (and can possibly
increase).
\end{corollary}

The normal form with non-increasing indices is obtained when we move
up in the weak Bruhat order by moving letters $x_i$ with ``small''
indices to the right of letters $x_j$ with ``large'' indices (we
assume $i<j$) by applying the substitution $x_ix_j \to x_{j+1}x_i$.
The normal form with unit decrease is obtained when we move down in
the weak Bruhat order by moving letters $x_j$ with ``very large''
indices to the right of letters $x_i$ with ``small'' indices (we
assume $j-1 > i$) by  applying the substitution $x_jx_i \to
x_ix_{j-1}$.

\begin{corollary}
The set of rules
\[ x_ix_j \to x_{j+1}x_i, \]
for $j<j$, represents a confluent rewriting system on $P_2$. The
same is true for the reversed set of rules
\[ x_{j+1}x_i \to x_ix_j, \]
for $i<j$.
\end{corollary}

We note here that both normal forms are well known, but the top one
is used more often in the literature on Thompson monoids and groups.
However, J.~Belk and K.~Brown use the bottom one quite efficiently
in~\cite{belk-b:forest,belk:phd} to get length functions for the
elements in $P_2$ and $F_2$ and then use these length functions in
further applications. Taking a different approach, J.~Belk shows
independently in~\cite{belk:phd} that that the two normal forms of
an element $g$ of $X$-length $n$ in $P_2$ bound the class of $X$
words representing $g$ in the so called word graph of $g$ based on
the rewriting rules above (without describing these classes as
unions of Tamari lattice congruence classes).

We observe that an inversion sequence word $x$ in $X_\infty$
represents a basis element of the free monoid $X_\infty/{\sim}$ from
Proposition~\ref{free} if and only if $x$ has a single occurrence of
$x_0$ (necessarily at the very end). The corresponding basis
elements in $S_\infty/{\sim}$ are represented by those permutations
that start with their largest term. One can easily pick normal
representatives (either top or bottom) for basis elements either in
$X_\infty/{\sim}$ or in $S_\infty/{\sim}$.

Finally, we mention that the linearizations corresponding to the
bottoms of Tamari classes on forests in $T^*$ are the ones obtained
by post-order linearizations, while the ones corresponding to the
tops are the inverse post-order linearizations.

\section{Polygon partitions and Thompson monoids}\label{partitions}

In this section we briefly indicate how, for $k \geq 2$, Tamari
orders on partitions of $((k-1)n+2)$-gons into $(k+1)$-gons lead
to Thompson monoid
\begin{equation}\label{presentation-pk}
 P_k = Mon \langle \ x_0,x_1,x_2,\dots \mid x_i x_j = x_{j+k-1}x_i, \text{ for }i<j \ \rangle.
\end{equation}

We fix $k$ and $n$ to discuss the general case, but provide
concrete examples with $k=3$ and $n=4$.

\begin{figure}[!ht]
\begin{center}
\includegraphics{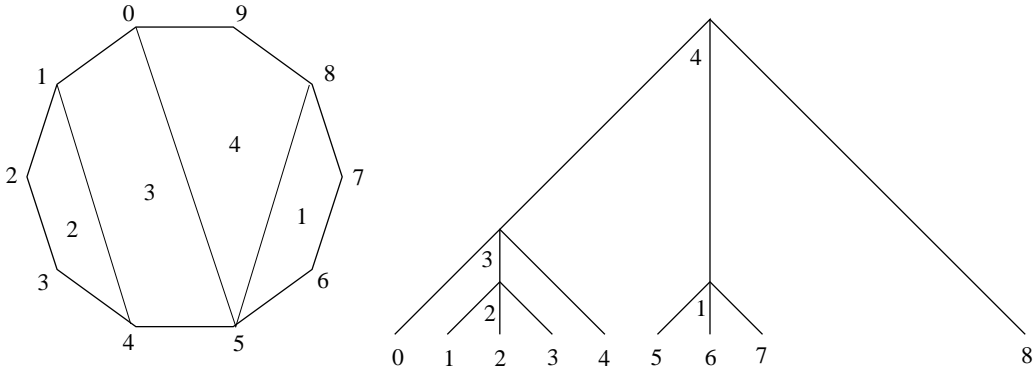}
\caption{A partition and a corresponding tree}\label{partition}
\end{center}
\end{figure}

First we define Tamari order on the partitions of a fixed
$((k-1)n+2)$-gon into $(k+1)$-gons. Label the vertices of the
$((k-1)n+2)$-gon by $0,\dots,(k-1)n+1$ in the positive direction.
Any diagonal $d$ used in the partition is common to two uniquely
determined $(k+1)$-gons in the partition that form a $2k$-gon using
$d$ as a diagonal (connecting opposite vertices in the $2k$-gon).
The labels of the $2k$-gon are still (cyclically) ordered from the
smallest to the largest in positive direction. Let the $k$ smallest
labels on the vertices of the $2k$-gon be $\ell_1,\dots,\ell_k$. The
diagonal $d$ has size $i$ if it uses the vertex labeled $\ell_i$. A
partition $Q_2$ covers a partition $Q_1$ if it is obtained from
$Q_1$ by removing a diagonal of size $i$ from a $2k$-gon in $Q_1$
and replacing it by the diagonal of size $i+1$ in the same $2k$-gon.
The Tamari partial order on partitions is then just the closure of
the cover relation. An example of  partition (with $k=3$, $n=4$) is
given in the left half of Figure~\ref{partition}. The diagonal
$(0,5)$ has size $1$, while the diagonals $(1,4)$ and $(5,8)$ have
size $2$.

Next we define linearized labeled $k$-ary rooted trees. A linearized
labeled rooted $k$-ary tree on $n$-interior vertices is a rooted
tree in which the root has degree $k$ (unless the tree has only the
root as a vertex, in which case its degree is 0), all interior
vertices have degree $k+1$, the leafs are labeled bijectively by $0,
\dots, (k-1)n$ and the interior vertices are labeled bijectively by
$0,\dots,n$ in such a way that the labels on each path from the root
to a leaf are decreasing. An example is given in the right half of
Figure~\ref{partition}. Denote the set of such trees by $T_{k,n}$.

We can define the notion of an interior vertex covering a gap just
as in the binary tree case. However, in this case every interior
vertex covers exactly $k-1$ gaps. For every tree $t$ in $T_{k,n}$
define a sequence $\pi(t)$  of length $(k-1)n$ by \eqref{pit}. In
our example in Figure~\ref{partition} we have $\pi(t)=32234114$.
Denote the image of $T_{k,n}$ by $S_{k,n}$. It consists of all
sequences of length $(k-1)n$ such that each term $1,2,3,\dots,n$
appears exactly $k-1$ times and all terms between two appearances of
a term $j$ are smaller than $j$. The map
\[ \pi: T_{k,n} \to S_{k,n} \]
is a bijection. The inverse map $\tau: S_{k,n} \to T_{k,n}$  can be
defined in a way analogous to the binary case. Namely, for $\sigma$
in $S_{k,n}$ the $k-1$ element set $\sigma^{-1}(i)$ will be called
the set of locations of the term $i$. Starting from $(k-1)n+1$
appropriately labeled leafs we add, in step $i$, a $k$-ary tree with
single interior vertex ($k$-caret) labeled by $i$ in such a way that
for each location $\ell$ of $i$ the gap $(\ell-1,\ell)$ is covered
by the interior vertex $i$.

From $S_{k,n}$ we can go by counting inversions to $X$-words of
length $n$. The only important remark is that the number of
inversions $\mathsf{inv}_i(\sigma)$ does not depend on the
particular occurrence of $i$ in $\sigma$. Thus we have a bijective
map $x: S_{k,n} \to X_{k,n}$, where $X_{k,n}$ is the set of words
\[
 X_{k,x} = \{ x_{i_1}x_{i_j} \dots x_{i_n} \mid 0 \leq i_j \leq (k-1)(n-j), \ j=1,\dots,n \ \}.
\]
Going back from $X_{k,n}$ to $S_{k,n}$ is accomplished by starting
from $(k-1)n$ open slots and then, in step $j$, placing $k-1$
copies of the term $j$ in consecutive available open slots after
leaving the first $i_j$ leftmost slots open.

The operation of concatenation of $X$-words still makes perfect
sense and leads to corresponding interlacing operation $\div$ on
sequences such as those in $S_{k,n}$ and stacking operation $\div$
on $k$-ary trees.

Every sequence in $S_{k,n}$ provides a way to build a polygon
partition. Start with the polygon with no edges or diagonals drawn.
Draw the edge $(0,(k-1)n+1$). In a sense that will be clear later
this is the root edge. Then, for $i=n,\dots,1$ (in that order!), in
step $i$ add the $k-1$ locations of $i$ in $\sigma$ to the path
(keep the vertices in the path always in increasing order). The
union of all the obtained paths is the desired partition.

In our running example $\sigma=32234114$ and we start with the
edge (0,9), then for $i=4$, we add the vertices 5 and 8, since
these are the locations of 4 in $\sigma$ and we obtain the path
$(0,5,8,9)$. Then for $i=3$, we add the vertices 1 and 4
(locations of 3) to get the path $(0,1,4,5,8,9)$. For $i=2$ we add
vertices 2 and 3 to get the path $(0,1,2,3,4,5,8,9)$ and finally
for $i=1$ we add 6 and 7 to get the path $(0,1,2,3,4,5,6,7,8,9)$.

For each $i=n\dots,1$, each time we add $k-1$ new vertices to the
path we add a new $(k+1)$-gon $K_i$ to the partition. If we keep the
label $i$ on $K_i$ we obtain linearized partitions. Dropping the
labels on $K_i$ amounts to de-linearization in the corresponding
trees and equivalence relation ${\sim}$ on $S_{k,n}$.

\begin{figure}[!ht]
\begin{center}
\includegraphics{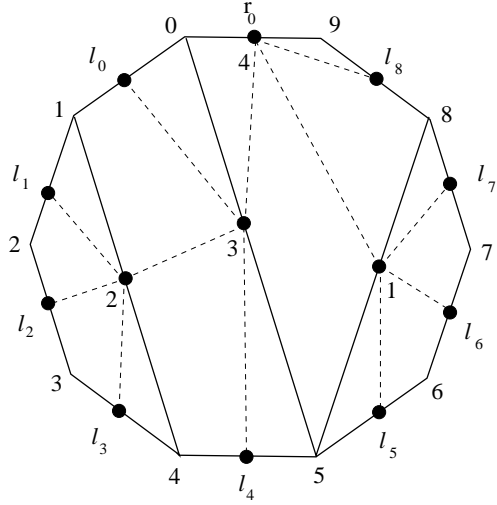}
\caption{The correspondence between trees and
partitions}\label{treepartition}
\end{center}
\end{figure}

There is a direct way to relate trees and partitions. Essentially,
the leafs $0,1,\dots,(k-1)n$ represent the edges
$(0,1),(1,2),\dots,((k-1)n,(k-1)n+1)$, the root represent the edge
$(0,(k-1)n+1)$ and the interior vertices $1,\dots,n-1$ represent the
diagonals. The interior vertex $n$, being the root, represents, the
edge $(0,(k-1)n+1)$. The diagram in Figure~\ref{treepartition}
depicts the correspondence in our model case for the tree and the
partition from Figure~\ref{partition}. The leafs are labeled by
$\ell_i$ and the root is labeled by $r_0$. The edges of the three
are dashed, while the partition edges are in full line. The vertices
of the tree are emphasized by representing them by small black
disks.

Note that the operation $\div$ on the level of polygon partitions
amounts to gluing partitioned polygons. Namely, if $Q_1$ is a
partition of an $(n_1+2)$-gon and $Q_2$ is a partition of an
$(n_2+2)$-gon into $(k+1)$-gons, then $R_1 \div R_2$ is obtained
by lifting all nonzero vertex labels in the second polygon by
$n_1$, and then gluing the two polygons along the edge $(0,n_1+1)$
in both polygons.

\begin{figure}[!ht]
\begin{center}
\includegraphics{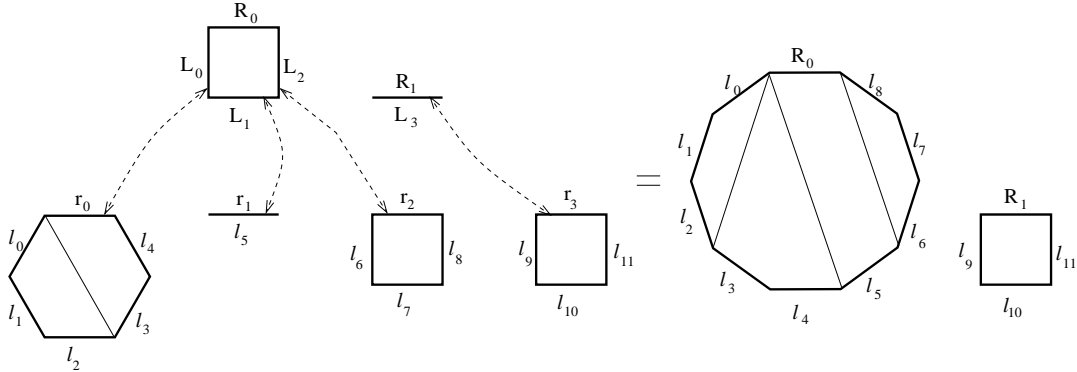}
\caption{Product of two partitions}\label{product-partition}
\end{center}
\end{figure}

Once again the situation can be lifted to arbitrary words in the
monoid $X^*$, which correspond to linearized $k$-ary forests in
$T^*$ with stacking operation, $*$-sequences with interlacing
operation, and partitions of finite sequences of polygons with
gluing operation. The operation on the level of polygon partitions
involves sequences of partitioned polygons and amounts to gluing the
root edges in the first sequence to the leaf edges with matching
label in the second partition. An example is given in
Figure~\ref{product-partition}.
The leaf edges in the first partition are labeled by $\ell_i$ and
the root edges by $r_i$, while capital letters are used in the
second partition. Note that trivial trees (single vertex, no edges)
in forests correspond to $2$-gons in polygons (represented as singe
edges). The dashed double arrows indicate which edges are to be
identified in the gluing process. One choice to describe the
multiplication of the depicted partitions by $*$-sequences is as
$2233**11*44 \div 11* = 22335511*44$. On the level of $X$-words the
above product corresponds to $x_6x_0x_0x_3 \div x_0 =
x_6x_0x_0x_3x_0$.

The Tamari equivalence ${\sim}$ is again a monoid congruence, the
top element in each class is always $132$-avoiding, the bottom one
is $231$-avoiding, and the obtained factor monoid is Thompson monoid
$P_k$, given by the presentation~\eqref{presentation-pk}. This is
just the positive monoid in the corresponding Thompson group $F_k$
given by the same presentation but as a group, which is the group
$F_k = P_k^{-1}P_k$ of fractions of $P_k$. The groups $F_k$, $k \geq
3$, share many properties with $F_2$
(see~\cite{brown:finiteness,brin-g:aut}).

Note that if $i$ is a concrete symbol that occurs on the left
(right) of some larger symbol $j$ (concrete or a $*$), then all
occurrences of $i$ are on the left (right) of $j$. Further, the
corresponding $k$-ary forest does not change when we exchange all
occurrences of $i$ and $i+1$ if and only if they are separated by
some larger symbol $j$ (thus $\overline{132}$ or $\overline{231}$
occurs - in the former case we go up and in the latter we go down
in the order).

The defining relations in $P_k$ have the form given
in~\eqref{presentation-pk} precisely because after all the
occurrences of $i$ and $i+1$ switch their places in a $*$-sequence
in which all $k-1$ occurrences of $i$ were to the left of all
$k-1$ occurrences of $i+1$ (separated by some larger symbol), the
number of inversions for the term $i$ increases by $k-1$.

\section{Concluding remarks}

It seems that Thompson monoid $P_2$ naturally codifies many
instances of Catalan-like objects in the sense that provides
``recipes'' for their construction as well as relations to
indicate which recipes lead to the same Catalan-like object.

The ubiquity of Catalan-like objects (enumerated by Catalan
numbers) is well known in combinatorics. On the other hand, the
ubiquity of Thompson's monoid (in fact the group) is equally well
known in infinite group theory. So it is fitting that these
objects are closely related. The fact that the finite Coxeter
groups of type $A$ play a role in the mix is also not extremely
surprising in the light of their own relevance in many situations.

It would be interesting to explore/establish connections between
Thompson monoids (not necessarily of type $F$) and Tamari lattices
of type $B$ and $D$~\cite{hugh:tamari-b,reading:cambrian}
(corresponding to factor lattices of finite Coxeter groups of type
$B$ and $D$).

The connection between the higher Thompson monoids $P_k$, $k \geq
3$, and the higher Catalan objects ($k$-ary forests) leads to a
natural question of exploring the sequences in $S_{n,k}$ as a kind
of higher Coxeter objects of type $A$ (ordinary permutations play
this role when $k=2$).

\section{Acknowledgments}
Thanks to Marcelo Aguiar, James Belk and Matt Brin for sharing their
thoughts. Also, thanks to Nathan Reading for very careful reading of
the text and numerous valuable suggestions and references.



\begin{thebibliography}{CFP96}

\bibitem[BB05]{belk-b:forest}
James Belk and Kenneth Brown.
\newblock Forest diagrams for elements of {T}hompson's group {$F$}.
\newblock {\em Internat. J. of Algebra Comput.}, 15(5-6):815--850, 2005.

\bibitem[Bel04]{belk:phd}
James Belk.
\newblock {\em Thompson's Group {$F$}}.
\newblock PhD thesis, Cornell University, 2004.

\bibitem[BG84]{brown-g:fp}
Kenneth~S. Brown and Ross Geoghegan.
\newblock An infinite-dimensional torsion-free {${\rm FP}\sb{\infty }$} group.
\newblock {\em Invent. Math.}, 77(2):367--381, 1984.

\bibitem[BG98]{brin-g:aut}
Matthew~G. Brin and Fernando Guzm{\'a}n.
\newblock Automorphisms of generalized {T}hompson groups.
\newblock {\em J. Algebra}, 203(1):285--348, 1998.

\bibitem[Bro87]{brown:finiteness}
Kenneth~S. Brown.
\newblock Finiteness properties of groups.
\newblock {\em J. Pure Appl. Algebra}, 44(1-3):45--75, 1987.

\bibitem[BS85]{brin-s:nofree}
Matthew~G. Brin and Craig~C. Squier.
\newblock Groups of piecewise linear homeomorphisms of the real line.
\newblock {\em Invent. Math.}, 79(3):485--498, 1985.

\bibitem[BW97]{bjorner-w:shellable2}
Anders Bj{\"o}rner and Michelle~L. Wachs.
\newblock Shellable nonpure complexes and posets. {II}.
\newblock {\em Trans. Amer. Math. Soc.}, 349(10):3945--3975, 1997.

\bibitem[CFP96]{cannon-f-p:thompson}
J.~W. Cannon, W.~J. Floyd, and W.~R. Parry.
\newblock Introductory notes on {R}ichard {T}hompson's groups.
\newblock {\em Enseign. Math. (2)}, 42(3-4):215--256, 1996.

\bibitem[ER96]{edelman-r:stasheff}
Paul~H. Edelman and Victor Reiner.
\newblock The higher {S}tasheff-{T}amari posets.
\newblock {\em Mathematika}, 43(1):127--154, 1996.

\bibitem[FH93]{freyd-h:hidempotents}
Peter Freyd and Alex Heller.
\newblock Splitting homotopy idempotents. {II}.
\newblock {\em J. Pure Appl. Algebra}, 89(1-2):93--106, 1993.

\bibitem[GS97]{guba-s:diagram}
Victor Guba and Mark Sapir.
\newblock Diagram groups.
\newblock {\em Mem. Amer. Math. Soc.}, 130(620):viii+117, 1997.

\bibitem[GS98]{grigorchuk-s:cancellation}
R.~I. Grigorchuk and A.~M. Stepin.
\newblock On the amenability of cancellation semigroups.
\newblock {\em Vestnik Moskov. Univ. Ser. I Mat. Mekh.}, (3):12--16, 73, 1998.

\bibitem[Gub06]{guba:quadratic}
Victor Guba.
\newblock The {D}ehn function of {R}ichard {T}hompson's group {F} is quadratic.
\newblock {\em Invent. Math.}, 163(2):313--342, 2006.

\bibitem[Hig74]{higman:simple}
Graham Higman.
\newblock {\em Finitely presented infinite simple groups}.
\newblock Department of Pure Mathematics, Department of Mathematics, I.A.S.
  Australian National University, Canberra, 1974.

\bibitem[HT72]{huang-t:lattice}
Samuel Huang and Dov Tamari.
\newblock Problems of associativity: {A} simple proof for the lattice property
  of systems ordered by a semi-associative law.
\newblock {\em J. Combinatorial Theory Ser. A}, 13:7--13, 1972.

\bibitem[Hug04]{hugh:tamari-b}
Thomas Hugh.
\newblock Tamari lattices and non-crossing partitions in types {$B$} and {$D$}.
\newblock Formal Power Series and Algebraic Combinatorics, 2004.

\bibitem[LR98]{loday-r:hopftrees}
Jean-Louis Loday and Mar{\'{\i}}a~O. Ronco.
\newblock Hopf algebra of the planar binary trees.
\newblock {\em Adv. Math.}, 139(2):293--309, 1998.

\bibitem[Rea]{reading:cambrian}
Nathan~P. Reading.
\newblock Cambrian lattices.
\newblock Adv.~Math., to appear.

\bibitem[Rea04]{reading:congruences}
Nathan~P. Reading.
\newblock Lattice congruences of the weak order.
\newblock {\em Order}, 21(4):315--344, 2004.

\bibitem[Sta97]{stanley:ec1}
Richard~P. Stanley.
\newblock {\em Enumerative combinatorics. {V}ol. 1}, volume~49 of {\em
  Cambridge Studies in Advanced Mathematics}.
\newblock Cambridge University Press, Cambridge, 1997.

\bibitem[Sta99]{stanley:ec2}
Richard~P. Stanley.
\newblock {\em Enumerative combinatorics. {V}ol. 2}.
\newblock Cambridge University Press, Cambridge, 1999.
\newblock With a foreword by Gian-Carlo Rota and appendix 1 by Sergey Fomin.

\bibitem[{\v{S}}un03]{sunik:sds}
Zoran {\v{S}}uni{\'k}.
\newblock Self-describing sequences and the {C}atalan family tree.
\newblock {\em Electron. J. Combin.}, 10:Note 5, 9 pp. (electronic), 2003.

\bibitem[Tho]{thompson:notes}
Richard~J. Thompson.
\newblock Handwritten notes, $\sim1960$'s.

\end{thebibliography}

\def\cprime{$'$}

\end{document}